\documentclass[leqno,12pt]{article}

\usepackage{latexsym}
 \usepackage{graphicx}
\usepackage{amsmath}
\usepackage{amssymb}
\usepackage{algorithm}
\usepackage{algorithmic}
\usepackage{cite}
\usepackage{color}

\usepackage[margin=1.3in, top=1.3in, bottom=1.3in]{geometry}

\newcommand{\nc}{\newcommand}
\nc{\nt}{\newtheorem}
\nt{thm}{Theorem}[section]
\nt{cor}[thm]{Corollary}
\nt{prop}[thm]{Proposition}
\nt{lem}[thm]{Lemma}
\nt{defn}[thm]{Definition}
\nt{rem}[thm]{Remark}
\nt{exa}[thm]{Example}
\nt{ass}[thm]{Assumption}
\nt{alg}[thm]{Algorithm}
\nc{\ip}[2]{\mbox{$\langle #1,#2 \rangle$}}
\nc{\pf}{\noindent{\bf Proof\ \ }}
\nc{\finpf}{\hfill{$\Box$}\linespace}
\nc{\linespace}{\vspace{\baselineskip} \noindent}
\nc{\R}{{\bf R}}
\nc{\bR}{{\overline{\R}}}
\nc{\E}{{\bf E}}
\nc{\U}{{\bf U}}
\nc{\V}{{\bf V}}
\nc{\W}{{\bf W}}
\nc{\X}{{\mathcal X}}
\nc{\Y}{{\mathcal Y}}
\nc{\M}{{\mathcal M}}
\nc{\Rn}{{\bf R}^n}
\nc{\bx}{\bar{x}}
\nc{\by}{\bar{y}}
\nc{\inT}{\mbox{\rm int}\,}
\nc{\cl}{\mbox{\rm cl}\,}
\nc{\gph}{\mbox{\rm gph}\,}

\def\tto{\;{\lower 1pt \hbox{$\rightarrow$}}\kern -12pt
           \hbox{\raise 2.8pt \hbox{$\rightarrow$}}\;}
\newenvironment{myequation}{\setcounter{equation}{\value{thm}}
   \begin{equation}}{\addtocounter{thm}{1}\end{equation}}

\nc{\bmye}{\begin{myequation}}
\nc{\emye}{\end{myequation}}

\begin{document}
\title{
Active-set Newton methods and \\
partial smoothness
}
\author{
A.S. Lewis
\thanks{ORIE, Cornell University, Ithaca, NY 14853, USA.
\texttt{people.orie.cornell.edu/aslewis} \hspace{2mm} \mbox{}
Research supported in part by National Science Foundation Grant DMS-1613996.}
\and
C. Wylie
\thanks{ORIE, Cornell University, Ithaca, NY 14853, USA.}
}
\date{\today}
\maketitle

\begin{abstract}
Diverse optimization algorithms correctly identify, in finite time, intrinsic constraints that must be active at optimality.  Analogous behavior extends beyond optimization to systems involving partly smooth operators, and in particular to variational inequalities over partly smooth sets.  As in classical nonlinear programming, such active-set structure underlies the design of accelerated local algorithms of Newton type.  We formalize this idea in broad generality as a simple linearization scheme for two intersecting manifolds.  
\end{abstract}
\medskip

\noindent{\bf Key words:} partial smoothness, active set identification, variational inequality
\medskip

\noindent{\bf AMS 2000 Subject Classification:} 90C31, 49M05, 65K10.
\medskip

\section{Introduction}
Active set methods for optimization and more general variational problems proceed, either explicitly or implicitly, in two phases.  The first phase relies on a globally convergent process to identify activity in the underlying constraints at the solution.  After correct identification, the problem simplifies locally:  the algorithm transitions to a second phase with accelerated local convergence, typically as a result of linearizing the simplified problem.  Standard texts like \cite{nocedal_wright,fac_pang} include overviews of active set methods for traditional nonlinear programming and variational inequalities, but the algorithmic philosophy also extends to less classical large-scale settings \cite{dual_av,prx_lin,liang-phd,fadili-activity}.

The idea of partial smoothness, introduced in \cite{Lewis-active}, formalizes the property driving the {\em first} phase:  ``identifiability''.  This term originated in \cite{Wright}, with precursors in \cite{Calamai-More87,Dunn87,Burke-More88,Burke90,Al-Khayyal-Kyparisis91,Ferris91,Flam92}.  Applications to  convergence and sensitivity analysis for large-scale convex optimization include \cite{fb_ps,activ_ident,liang,Vaiter2017,fadili-activity,fadili-model,mirror}.  Rather than active sets of  constraints, partly smooth geometry focuses on the ``active manifold'' defined by those constraints.  The equivalence with identifiability is explored in detail in \cite{ident}.

In contrast with most of this prior literature, in this work we concentrate on the {\em second} phase of active-set methods:  in partly smooth settings, how should we understand accelerated local convergence?  A formalism parallel to  partial smoothness with related algorithmic ambitions, ``$\mathcal{VU}$ theory'', emerged in \cite{LS97}, and includes \cite{LOS,MC02,MC05,MC04,HS09,MC03}, and most recently the survey \cite{claudia-icm}.  
Algorithms exploiting the active manifold structure, and the resulting connections between partial smoothness, $\mathcal{VU}$ theory, and classical nonlinear programming were explored in \cite{miller-malick}.
Our aims here, rooted in the partly smooth framework, are more general, formal and geometric.

In essence, as explained in \cite{lewis-liang}, partial smoothness is a constant-rank property.  Assuming this simple differential-geometric condition, the underlying problem reduces, locally, to finding the unique point in two manifolds that intersect transversally.  Just like the classical Newton method for solving equations, via linearization we arrive at fast local algorithms.

A standard illustration (see \cite{lewis-liang}) involves $C^{(2)}$-smooth functions $g_i \colon \Rn \to \R$, for indices $i=1,2,3\ldots,m$, and a corresponding set
\[
F ~=~ \big\{x : g_i(x) \le 0 ~ \forall i \big\}.
\]
Given a cost vector $\bar c \in \Rn$, the first-order optimality condition for the optimization problem 
$\sup_{x \in F} \ip{\bar c}{x}$ is a generalized equation involving the normal cone to $F$:
\bmye \label{optimality}
\bar c \in N_F(x).
\emye
Consider any feasible point $x \in F$ for which the constraint gradients $\nabla g_i(x)$ for $i$ in the active set
\[
I(x) ~=~ \{i : g_i(x) = 0\}
\]
are linearly independent:  for such points, we can describe the normal cone operator using Lagrange multiplier vectors $\lambda \in \R^m$:
\[
N_F(x) ~=~ \Big\{ \sum_i \lambda_i \nabla g_i(x) : \lambda \ge 0,~ \lambda_i = 0 ~\big(i \not\in I(x)\big) \Big\}.
\]
Let us focus on a particular such point, $\bar x$, that solves the equation (\ref{optimality}), with corresponding multiplier $\bar\lambda$.  Assuming strict complementarity (meaning $\bar\lambda_i > 0$ for all $i \in I(\bar x)$), the operator $N_F$ is partly smooth at $\bar x$ (for $\bar c$):  around the pair 
$(\bar x , \bar c)$, its graph, $\gph N_F$, is a manifold, and the canonical projection is constant-rank and maps $\gph N_F$ locally onto the active manifold
\[
{\mathcal M} ~=~ \Big\{ x : g_i(x) = 0 ~ \big(i \in I(\bar x) \big) \Big\}.
\]
In particular, for any nearby pair $(x,c)$ in the graph, the point $x$ must lie in ${\mathcal M}$, so has the same active set:  this is the property of identifiability.  For more details in the setting of this standard example, see \cite{lewis-liang}.

Now, however, suppose in addition that the given point $\bar x$ also satisfies the standard {\em second-order} sufficient conditions:  the Lagrangian $\sum_i \bar\lambda_i \nabla g_i - c$ is not only critical at $\bar x$, but also its Hessian is positive definite on the tangent space to ${\mathcal M}$:
\[
\Big\{ d \in \Rn : \ip{\nabla g_i(\bar x)}{d} = 0 ~\big(i \in I(\bar x)\big) \Big\}.
\]
In that case, locally around the pair $(\bar x, \bar c)$, the graph of the normal cone operator, 
$\gph N_F$, and the space $\Rn \times \{\bar c\}$ are manifolds intersecting at the unique point 
$(\bar x, \bar c)$, and furthermore the intersection is transversal.  Based at a nearby point 
$(x,c) \in \gph N_F$, we could linearize $\gph N_F$, calculate the intersection point of the resulting tangent space with $\Rn \times \{\bar c\}$, restore this point to $\gph N_F$, and repeat.  As we discuss later, such a process fundamentally underlies the philosophy of classical sequential quadratic programming.

In this classical example, we see that our theme just results in a standard methodology, namely sequential quadratic programming for nonlinear programming.  Nonetheless, as should be apparent, even the concise summary we gave inevitably involves an extended set-up involving active constraint indices, Lagrange multipliers, strict complementarity, and second-order sufficient conditions.  Rather than algorithmic innovation, in this case our contribution is conceptual.    By contrast with the traditional development, our theme is simpler and geometric, involving just fundamental ideas about manifolds:  tangent spaces, transversality, and constant-rank maps.  In summary, while we hope that the ideas we develop here may help illuminate numerical methods in less traditional settings, our primary goal is not the immediate design of new algorithms.  Rather our aim is a simple unifying framework in which to study active-set-type acceleration techniques for a much broader class of optimization and variational problems.

\section{Manifolds}
We begin our more formal development with a brief summary of terminology for manifolds.  We refer the reader to \cite{Lee} for basic background in differential geometry.

We consider manifolds in a Euclidean space $\E$, by which we mean sets defined, locally, by smooth equations with linearly independent gradients.  More precisely, a set 
$\X \subset \E$ is a $C^{(r)}$-{\em manifold} around a point 
$z \in \X$ (for some degree $r = 1,2,3,\ldots$) if there exists a $C^{(r)}$-smooth map $F \colon \E \to \R^m$ (for some number $m$) with surjective derivative $DF(z)$ such that the inverse image $F^{-1}(0)$ is a neighborhood of $z$ in $\X$.  Equivalently, there exists a number $k$, an open neighborhood $U$ of zero in $\R^k$, and a $C^{(r)}$-smooth map $G \colon U \to \E$ with value $G(0) = z$, everywhere injective derivative, and range $G(U)$ a neighborhood of $z$ in $\X$.  In that case the number $k = \dim \E - m$ is called the {\em dimension} of the the manifold $\X$ around $z$, and we refer to $G$ as a {\em local coordinate system}. 

With the notation above, the {\em tangent} and {\em normal} spaces to $\X$ at $z$ are given respectively by
\begin{eqnarray*}
T_\X(x) &=& \mbox{Null}\big(DF(z)\big) ~=~ \mbox{Range}\big(DG(0)\big) \\
N_\X(x) &=& \mbox{Range}\big(DF(z)^*\big) ~=~ \mbox{Null}\big(DG(0)^*\big). 
\end{eqnarray*}
These two orthogonal complements are independent of the choice of maps $F$ and $G$.  We say that two $C^{(2)}$-manifolds around $z$ intersect {\em transversally} there if the two corresponding normal spaces  intersect trivially.

\section{A Newton method for intersecting manifolds}
The following simple scheme uses standard ideas.  We study the intersection of two manifolds by linearizing just one of them.

\begin{thm}[Semi-linearization error]
\label{quadratic}
In a Euclidean space, consider an isolated point of intersection $z$ of two \mbox{$C^{(2)}$-manifolds} 
$\X$ and $\Y$ around $z$, and suppose that the intersection is transversal.  Then there exists an open neighborhood $W$ of $z$ such that, for all points $x \in \X \cap W$, the linear approximation $x + T_\X(x)$ intersects the set $\Y \cap W$ at a unique point $y(x)$, which satisfies
\bmye \label{small-order}
|y(x)-z| = O(|x-z|^2) \quad \mbox{as}~ x \to z.
\emye
\end{thm}

\pf
Since the manifolds $\X$ and $\Y$ intersect at a singleton, and transversally, their dimensions sum to $\dim\E$.  We can represent $\X$ using a local coordinate system:  an open neighborhood $U$ of zero in $\R^k$ and a $C^{(2)}$-smooth map $G \colon U \to \E$ with value $G(0) = z$, everywhere injective derivative, and range $G(U)$ an open neighborhood of $z$ in $\X$.  By definition, there also exists a $C^{(2)}$-smooth map 
$F \colon \E \to \R^k$ with surjective derivative $DF(z)$ such that the inverse image $F^{-1}(0)$ is a neighborhood of $z$ in $\Y$.

The transversality condition is
\[
\{0\} ~=~ N_\X(z) \cap N_\Y(z) ~=~ \mbox{Null}\big(DG(0)^*\big) \cap \mbox{Range}\big(DF(z)^*\big).
\]
Hence all vectors $v \in \R^k$ satisfy
\[
DG(0)^* \big(DF(z)^* v\big) = 0 \quad \Leftrightarrow \quad DF(z)^* v = 0.
\]
By assumption the adjoint map $DF(z)^*$ is injective, so the right-hand side is equivalent to $v= 0$.  Consequently the composition $DG(0)^* \circ DF(z)^*$ is injective, so its adjoint
$DF(z) \circ DG(0)  \colon \R^k \to \R^k$ is surjective and hence in fact invertible.  By the inverse function theorem, the composite map $F \circ G \colon \R^k \to \R^k$ has a $C^{(2)}$-smooth inverse map around zero.

We can express any point $x \in \X$ near $z$ as $x=G(u)$ for some small vector $u \in \R^k$.  In that case, the point we seek has the form $y(x) = G(u) + DG(u)v$ for some small vector $v \in \R^k$, where
\bmye \label{intersection}
F\big(G(u) + DG(u)v\big) = 0.
\emye
The function of $(u,v)$ on the left-hand side, which we denote $H(u,v)$, is $C^{(1)}$-smooth.   Its value at $(u,v)=(0,0)$ is zero, and its derivative there with respect to $v$ is 
$D_v H(0,0) = DF(z) \circ DG(0)$.  Since this derivative is invertible, the implicit function theorem implies that, for all small vectors $u \in \R^k$, equation (\ref{intersection}) has a unique small solution $v(u) \in \R^k$, depending smoothly on $u$, and satisfying $v(0) = 0$.

Since the map $G$ is $C^{(2)}$-smooth, we have
\bmye \label{expansion}
G(u+v) = G(u) + DG(u)v + O(|v|^2)
\emye
for all small vectors $u$ and $v$.  Since the map $F$ is Lipschitz, we deduce
\begin{eqnarray*}
F\Big(G\big(u+v(u)\big)\Big) 
&=& 
F\Big( G(u) + DG(u)v(u) + O\big(|v(u)|^2\big) \Big) \\
&=& 
F\Big( G(u) + DG(u)v(u)\Big) + O\big(|v(u)|^2\big) \\
&=& 
O\big(|v(u)|^2\big),
\end{eqnarray*}
as $u \to 0$, and hence
\[
u+v(u) = O\big(|v(u)|^2\big),
\]
since the inverse map $(F \circ G)^{-1}$ is Lipschitz.  But since
\[
|v(u)| \le |u| + |u+v(u)| = |u| + O\big(|v(u)|^2\big) \le |u| + \frac{1}{2}|v(u)|
\]
for all small $u$, we deduce $|v(u)| \le 2|u|$.  Equation (\ref{expansion}) implies
\begin{eqnarray*}
y(x) - z &=& G(u) + DG(u)v(u) - z ~=~ G\big(u+v(u)\big) + O\big(|v(u)|^2\big) - z  \\
&=& G\big(O\big(|v(u)|^2\big) - z + O\big(|v(u)|^2\big) \\
&=& O\big(|v(u)|^2\big) 
~=~ O(|u|^2) 
~=~ O\big(|G(u)-z|^2\big)
~=~ O(|x-z|^2)
\end{eqnarray*}
as claimed.
\finpf

If we can somehow restore the point $y(x)$ in the manifold $\Y$ to a nearby point in the manifold $\X$, by iterating we rapidly converge to the intersection point $z$.  Based on this idea, the following simple result is the foundation for all the Newton schemes we consider.

\begin{cor}[Semi-linearization for intersecting manifolds] \label{semi}
Following the framework of Theorem \ref{quadratic}, consider two manifolds $\X$ and $\Y$ intersecting transversally at their unique common point $z$, and in addition a Lipschitz map 
$R \colon \Y \to \X$ satisfying $R(z)=z$.  Then, starting from any point 
$x \in \X$ near $z$, the iteration
\[
x \leftarrow R\big(y(x)\big)
\]
converges quadratically to $z$.
\end{cor}

\pf
Using the fixed point and Lipschitz properties, and equation (\ref{small-order}), we deduce
\[
\big|R\big(y(x)\big) - z\big| = \big|R\big(y(x)\big) - R(z)\big| = O(|y(x)-z|) = O(|x-z|^2),
\]
for all points $x \in \X$ near $z$, as required.
\finpf

Since $\X$ is a $C^{(2)}$-smooth manifold around the point $z$, a natural candidate for the restoration map $R$ is $\mbox{Proj}_\X$, the corresponding nearest-point projection map onto $\X$, which is $C^{(1)}$-smooth and hence Lipschitz around $z$, and which obviously leaves $z$ fixed.  However, in the broad scenario that we have in mind, this projection may not be the most natural choice.  We describe this scenario next.

\section{Generalized equations} \label{generalized-equations}
Henceforth we fix a Euclidean space $\U$, and consider the problem of finding a zero of a set-valued mapping $\Phi$, or in other words a solution of the generalized equation
\[
0 \in \Phi(u).
\]
Henceforth we rely on standard terminology from variational analysis, as presented in monographs such as \cite{VA,Mord_1}.

Assuming that the graph of $\Phi$ is a manifold, locally, we can frame this equation as an intersection problem as in the previous section and Theorem \ref{quadratic}.  The underlying Euclidean space is now $\E = \U^2$, the sets $\X$ and $\Y$ are the graphs of $\Phi$ and the zero mapping respectively,
\begin{eqnarray*}
\X &=& \mbox{gph}\,\Phi ~=~ \{(u,v) \in \U^2 : v \in \Phi(u) \} \\
\Y &=& \U \times \{0\},
\end{eqnarray*}
and the intersection point is $z=(\bar u,0)$.  The transversality condition is
\[
N_{\mbox{\scriptsize gph}\,\Phi}(\bar u,0) \cap N_{\U \times \{0\}}(\bar u, 0) = \{(0,0)\}.
\]
Automatically, then, we must have $\dim \mbox{gph}\,\Phi = \dim\U$.  

Using the language of {\em coderivatives}, we have by definition for any set-valued mapping $\Phi$ the relation
\[
(v, -y) \in N_{\mbox{\scriptsize gph}\,\Phi}(\bar u,0) \quad \Leftrightarrow \quad
v \in D^*\Phi (\bar u | 0)(y).
\]
In this definition, $N$ denotes the {\em limiting normal cone}, which coincides with the classical normal space for manifolds.  Hence, since
\[
N_{\U \times \{0\}}(\bar u, 0) = \{0\} \times \U,
\]
we can rewrite the transversality condition as in the following assumption.

\begin{ass}
{\rm
The point $\bar u \in \U$ is an isolated zero of a set-valued mapping $\Phi \colon \U \tto \U$.  Furthermore, the graph of $\Phi$ is a $C^{(2)}$-smooth manifold around the point $(\bar u,0)$, and the transversality condition
\bmye \label{transversality}
0 \in D^*\Phi (\bar u | 0)(y) \quad \Rightarrow \quad y=0.
\emye
holds.
}
\end{ass}

\noindent
We remark in passing that the coderivative criterion (\ref{transversality}) plays a fundamental role throughout the variational analysis literature.  For any set-valued mapping $\Phi$ with closed graph and satisfying $0 \in \Phi(\bar u)$, it is equivalent to {\em metric regularity} \cite[Theorem 9.43]{VA}:  the ratio of two distances
\[
\frac{\mbox{dist}\big(u,\Phi^{-1}(v)\big)}{\mbox{dist}\big(v,\Phi(u)\big)}
\]
is uniformly bounded above for pairs $(u,v)$ near $(\bar u, 0)$ (where we interpret $\frac{0}{0} = 0$).

Turning to the linearization step in Theorem \ref{quadratic}, given a current point $u \in \U$ and a current value $v \in \Phi(u)$, using the language of {\em graphical derivatives}, we have by definition 
\bmye \label{Newton-step}
(u',0) \in (u,v) + T_{\mbox{\scriptsize gph}\,\Phi}(u,v)
\quad \Leftrightarrow \quad 
-v \in D\Phi (u | v)(u'-u).
\emye
In this definition, $T$ denotes the standard {\em tangent cone}, which coincides with the classical tangent space for manifolds.  Having calculated the new point $u' \in \U$ from this relationship, to apply Corollary \ref{semi} we need a Lipschitz map to restore the point $(u',0)$ to a nearby point on 
$\gph\,\Phi$.

\subsubsection*{Example:  classical Newton method}
For illustration, consider the case where the mapping $\Phi$ is just a single-valued $C^{(2)}$-smooth map 
$F \colon \U \to \U$.  The transversality condition (\ref{transversality}) reduces to invertibility of the derivative map $DF(\bar u)$.  At the current pair $(u,F(u)) \in \mbox{gph}\,F$, the linearization step (\ref{Newton-step}) results in a new pair $(u',0) \in \U \times \{0\}$ where 
$u' = u -( DF(u))^{-1}F(u)$.  Rather than trying to project $(u',0)$ onto $\mbox{gph}\,F$, we can use the natural restoration map $R$ defined by $(u',0) \mapsto \big(u',F(u')\big)$, which fixes the point 
$(\bar u,0)$ and obviously inherits the Lipschitz property of $F$.  Corollary \ref{semi} now simply recovers the quadratic convergence of the classical Newton iteration for solving $F(u)=0$.

\subsubsection*{Example:  variational inequalities}
More generally, consider the case
\[
\Phi = F + \Psi,
\]
for a single-valued $C^{(2)}$-smooth map $F \colon \U \to \U$ and a maximal monotone operator 
$\Psi \colon \U \tto \U$.  In particular, $\Psi$ might be the subdifferential operator $\partial g$ for a proper closed convex function $g \colon \U \to \overline\R$.  This case includes two important examples:  the special case when $F$ is the derivative $\nabla f$ of a smooth function $f \colon \U \to \R$, corresponding to the optimization problem $\inf\{f+g\}$, and the special case when
$g$ is the indicator function for a closed convex set 
$K \subset \U$, corresponding to the {\em variational inequality}
\[
\mbox{VI}(K,F): \quad 0 \in F(u) + N_K(u).
\]

We assume that the point $\bar u \in \U$ is an isolated zero of the mapping $F+\Psi$, and that the graph of the operator $\Psi$ is a $C^{(2)}$-smooth manifold around the point $\big(\bar u,-F(\bar u)\big)$.  Standard nonsmooth calculus (or direct calculation) implies
\begin{eqnarray*}
D^*(F+\Psi)(\bar u | 0)(y) &=& DF(\bar u)^*y + D^*\Psi\big(\bar u | -F(\bar u)\big)(y) \\
D(F+\Psi)(u | F(u) + w)(s) &=& DF(u)s + D\Psi\big(u | w\big)(s)
\end{eqnarray*}
for $w \in \Psi(u)$.
Hence the transversality condition (\ref{transversality}) becomes
\[
-DF(\bar u)^*y \in D^*\Psi \big(\bar u | -F(\bar u)\big)(y) \quad \Rightarrow \quad y=0.
\]
Given a current vector $u \in \U$ and a vector $w \in \Psi(u)$ corresponding to the choice $v = F(u)+w$ in  the relation (\ref{Newton-step}), the linearization step is defined by
\[
-F(u) - w ~\in~ DF(u)(u'-u) + D\Psi\big(u | w)\big)(u'-u).
\]

It remains to find a Lipschitz map restoring the new point $(u',0)$ to $\mbox{gph}\,(F+\Psi)$.  To that end we can use Minty's Theorem \cite[Theorem 12.15]{VA}, which guarantees that, for any choice of 
$\lambda > 0$, the {\em resolvent} map 
$(I + \lambda \Psi)^{-1} \colon \U \to \U$ is single-valued and nonexpansive.  Define a Lipschitz map
$Q \colon \U \to \U$ by
\bmye \label{restoration-Q}
Q(u) ~=~ (I + \lambda \Psi)^{-1}\big(u - \lambda F(u)\big).
\emye
Clearly zeros of the mapping $F+\Psi$ coincide with fixed points of $Q$, and
simple calculations show 
\[
\Big(Q(u) , \frac{1}{\lambda}\big(u  - Q(u)\big) - F(u) + F\big(Q(u)\big) \Big) 
~\in~ \gph(F+\Psi) \quad \mbox{for all}~ u \in \U.
\]
Hence we can use the restoration map $R \colon \U \times \{0\} \to \gph(F+\Psi)$ defined by setting 
$R(u,0)$ to the left-hand side of this inclusion.  This map is clearly Lipschitz, and fixes the point 
$(\bar u,0)$, as required.

In the smooth case ($\Psi=0$), we arrive at a slight variant of the classical Newton method for the smooth equation $F(u)=0$ that we saw in the previous example (and coinciding with it in the limiting case 
$\lambda = 0$).  A simple extension is when the mapping $\Psi$ is the normal cone operator $N_L$ for a linear subspace $L \subset \U$, in which case we arrive at a slight variant of the natural projected Newton scheme for solving $F(u) = 0,~ u \in L$.

\section{Identifying primal activity:  partial smoothness}
Corollary \ref{semi} describes the basic Newton-type approach to the manifold intersection scenario
\[
\X \cap \Y = \{z\}
\]
and more particularly for the special case associated generalized equation $0 \in \Phi(u)$:
\[
\mbox{gph}\,\Phi \cap (\U \times \{0\}) ~=~ \{(\bar u,0)\}.
\]
The Newton iteration depends fundamentally on a Lipschitz restoration map $R$ from the manifold $\X$ into the manifold $\Y$ that fixes the intersection point $z$.  In general it may be unclear how to construct such a map, and even in special cases when a closed form is available, as in the previous maximal monotone examples, the resulting conceptual algorithm is not transparent, and involves simultaneous coupled updates to primal variables (in the domain of the mapping $\Phi$) and dual variables (in the range of $\Phi$).

In this section we describe a systematic approach to the restoration map $R$, motivated by traditional active-set optimization algorithms.  We decouple updates to the primal and dual variables, thereby decomposing $R$ into a composition of two maps.  The first map updates the primal variables by identifying an underlying active submanifold in the primal manifold $\Y = \U \times \{0\}$.  The second map fixes the primal variables and updates the dual variables to restore the iterate to the manifold $\X = \mbox{gph}\,\Phi$.  This approach depends crucially on the existence of a constant-rank map from $\X$ into $\Y$:  in the special case of the generalized equation $0 \in \Phi(u)$, this assumption amounts to the ``partial smoothness'' of the mapping $\Phi$, language that we argue elsewhere is convenient for the philosophy of active set identification.

\begin{thm}[Active set method for manifold intersection] \label{active-intersection}
\label{newton}
In a Euclidean space, consider an isolated point of intersection $z$ of two \mbox{$C^{(2)}$-manifolds} 
$\X$ and $\Y$ around $z$, and suppose that the intersection is transversal.  Suppose furthermore that there exists a map $P \colon \X \to \Y$ that fixes the point $z$ and is $C^{(2)}$-smooth and constant rank nearby.  Then, for any sufficiently small neighborhood $W$ of $z$, the following properties hold:
\begin{enumerate}
\item[{\rm (a)}]
The image $\M = P(\X \cap W) \subset \Y$ is a $C^{(2)}$-manifold around $z$.
\item[{\rm (b)}]
There exists a Lipschitz map $S \colon \Y \cap W \to \M$ that fixes the point $z$.  In particular, we could choose $S$ to be the projection map
\[
\mbox{\rm Proj}_\M \colon \Y \cap W \to \M.
\]
\item[{\rm (c)}]
There exists a Lipschitz map $Q \colon \M \cap W \to \X$ that fixes the point $z$ and is a right inverse for $P$:
\[
P\big(Q(w)\big) = w \quad \mbox{for all}~ w \in \M \cap W.
\]
In particular, we could define
\[
Q(w) = \mbox{\rm Proj}_{P^{-1}(w)}(w) \quad \mbox{for all}~ w \in \M \cap W.
\]
\item[{\rm (d)}]
For all points $x \in \X \cap W$, the linear approximation $x + T_\X(x)$ intersects the set $\Y \cap W$ at a unique point $y(x)$.
\end{enumerate}
Furthermore, starting from any point $x \in \X$ near $z$, the iteration
\[
x \leftarrow Q\Big(S\big(y(x)\big)\Big)
\]
converges quadratically to $z$.
\end{thm}

\pf
Property (a) follows from the Constant Rank Theorem.  Property (b) follows from the well-known fact that the projection onto a $C^{(2)}$-manifold around a point $z$ is $C^{(1)}$-smooth around $z$.  Property (d) and the final conclusion follow from Corollary~\ref{semi}.  To complete the proof, we show that the map 
$Q$ defined in part (c) is in fact $C^{(1)}$-smooth around $z$.

By the Constant Rank Theorem, there exist local coordinate systems for the manifolds $\X$ and $\Y$ around the point $z$ with respect to which the map $P$ is linear, and in fact maps the point with coordinate vector $(u,v) \in \U \times \V$ to the point with coordinate vector $(u,0) \in  \U \times \W$, for suitable Euclidean spaces $\U$, $\V$, and $\W$.  To be more precise, there exist \mbox{$C^{(2)}$-smooth} maps 
$G \colon \U \times \V \to \E$ and $H \colon \U \times \W \to \E$ (the coordinate systems for $\X$ and $\Y$ respectively), satisfying $G(0,0) = z = H(0,0)$ and with injective derivatives $DG(0,0)$ and $DH(0,0)$, and such that
\[
P\big(G(u,v)\big) = H(u,0) \quad \mbox{for all small}~ u \in \U,~ v \in \V.
\]
Around $z$, the manifold $\M$ then just consists of points of the form $H(u,0)$ for small vectors 
$u \in \U$.  An obvious example of a right inverse map $Q$ with the desired properties is given by
\[
Q\big(H(u,0)\big) = G(u,0)  \quad \mbox{for all small}~ u \in \U,
\]
since then we have
\[
P\Big(Q\big(H(u,0)\big)\Big) = P\big(G(u,0)\big) = H(u,0)  \quad \mbox{for all small}~ u \in \U.
\]

For any small $\delta > 0$, the coordinate map $G$ is a diffeomorphism from the compact set
\[
\{(u,v) \in \U \times \V : |u| \le \delta,~ |v| \le \delta \}
\]
onto a compact neighborhood $\hat\X$ of the point $z$ in the manifold $\X$, and the map $u \mapsto H(u,0)$ is a diffeomorphism from the the compact set $\{u \in \U : |u| \le \delta\}$ onto a compact neighborhood of $z$ in the manifold $\M$.  Furthermore, for all $u \in \U$ satisfying $|u| \le \delta$ we have
\[
P^{-1}\big(H(u,0)\big) \cap \hat\X ~=~ \{G(u,v) : v \in \V,~ |v| \le \delta \}.
\]
Thus for all small $u$, projecting the point $w = H(u,0)$ onto the set $P^{-1}(w)$ amounts to solving the optimization problem
\[
\inf_{v \in \V} \big\{ |G(u,v) - H(u,0)|^2 : |v| \le \delta \big\}.
\]
The infimum is attained, by compactness.  By an easy continuity argument, if $u$ is sufficiently small, any optimal $v$ satisfies $|v| < \delta$, and hence also the first-order necessary condition
\bmye \label{necessary}
D_vG(u,v)^*\big(G(u,v) - H(u,0)\big) = (0,0).
\emye
The left-hand side of this equation is $C^{(1)}$-smooth.  When $u=0$, the equation reduces to
\[
D_vG(0,v)^*\big(G(0,v)\big) = (0,0).
\]
The derivative at $v=0$ of the map on the left-hand side is $D_vG(0,0)^*D_vG(0,0)$, which is invertible since $D_vG(0,0)$ is injective.  By the implicit function theorem, for all small $u \in \U$, equation (\ref{necessary}) has a unique small solution $v(u)$, depending smoothly on $u$, and clearly $v(0)=0$.
Since we have shown, for such $u$, the property
\[
\mbox{\rm Proj}_{P^{-1}(H(u,0))}\big(H(u,0)\big) = G\big(u,v(u)\big), 
\]
the result now follows.
\finpf

Returning to the setting of generalized equations $0 \in \Phi(u)$, for a variable $u$ in a Euclidean space $\U$ and a set-valued mapping $\Phi \colon \U \tto \U$, we next introduce the notion of partial smoothness  developed in \cite{lewis-liang}.

\begin{defn}[Partly smooth mappings] \label{partly}
{\rm
For any degree $r = 1,2,3,\ldots$, a set-valued mapping $\Phi \colon \U \tto \U$ is called $C^{(r)}$-{\em partly smooth at} a point 
$\bar u \in \U$ {\em for} a value $\bar v \in \Phi(\bar u)$ when its graph is a $C^{(r)}$-smooth manifold around the point
$(\bar u, \bar v)$ and the map defined by $(u,v) \mapsto u$ for $(u,v) \in \gph\Phi$ is constant rank nearby.  
}
\end{defn}

\noindent
By definition, the constant rank condition means that the coderivative
\[
D^*\Phi(u,v)(0) ~=~ \big\{ w \in \U : (w,0) \in N_{\mbox{\scriptsize \gph}\,\Phi}(u,v) \big\},
\]
is a subspace of constant dimension for points $(u,v)$ near $(\bar u,\bar v)$.  The following result is \cite[Proposition 3.2]{lewis-liang}.

\begin{prop}[Active manifold]
Assuming Definition \ref{partly}, there exists a set $\M$ in the space $\U$ such that any point $(u,v)$ sufficiently near the point $(\bar u,\bar v)$ in the graph of the mapping $\Phi$ satisfies $u \in \M$.  Any such set $\M$ is a $C^{(r)}$-smooth manifold around $\bar u$ and is locally unique:  that is, any two such sets are identical around $\bar u$.
\end{prop}

\noindent
We call the set $\M$ in this result the {\em active manifold at} $\bar u$ {\em for} $\bar v$.

With this terminology, we can specialize Theorem \ref{active-intersection} to generalized equations.

\begin{thm}[Active set method for generalized equations] \label{active-set} \hfill \mbox{}
Consider an isolated zero $\bar u \in \U$ for a set-valued mapping 
$\Phi \colon \U \tto \U$.  Suppose that $\Phi$ is $C^{(2)}$-partly smooth at $\bar u$ for $0$, with active manifold $\M \subset \U$, and assume the transversality condition
\[
0 \in D^*\Phi (\bar u | 0)(y) \quad \Rightarrow \quad y=0.
\]
Then, for any sufficiently small neighborhood $W$ of $\bar u$, the following properties hold:
\begin{itemize}
\item
There exists a Lipschitz map $S \colon W \to \M$ that fixes the point $\bar u$.  In particular, we could choose $S$ to be the projection map
\[
\mbox{\rm Proj}_\M \colon W \to \M.
\]
\item
There exists a Lipschitz map $Q \colon \M \cap W \to \U$ satisfying $Q(\bar u)=0$ and
\[
Q(u) \in \Phi(u) \quad \mbox{for all}~ u \in \M \cap W.
\]
In particular, we could define
\[
Q(u) = \mbox{\rm Proj}_{\Phi(u)}(0) \quad \mbox{for all}~ u \in \M \cap W.
\]
\item
For all points $(u,v) \in \gph\Phi$ near $(\bar u,0)$, the linearized equation
\[
-v \in D\Phi (u | v)(w-u), \quad w \in W
\]
has a unique solution $w(u,v)$.
\end{itemize}
Furthermore, starting from any point $(u,v) \in \gph \Phi$ near $(\bar u,0)$, the iteration
\[
u \leftarrow S\big(w(u,v)\big), \quad v \leftarrow Q(u)
\]
converges quadratically to $(\bar u,0)$.
\end{thm}

\pf
We simply apply Theorem \ref{newton} in the framework of Section \ref{generalized-equations}:
\begin{eqnarray*}
\X = \gph\Phi, & &  \Y = \U \times \{0\} \\
P(u,v)=(u,0) & \mbox{for} & (u,v) \in \gph\Phi.
\end{eqnarray*}
The result follows immediately.
\finpf

\section{Partly smooth variational inequalities}
Consider an isolated solution $\bar u$ of the variational inequality 
$\mbox{VI}(K,F)$:
\[
0 \in F(u) + N_K(u),
\]
for a $C^{(2)}$-smooth map $F \colon \U \to \U$ and a closed (possibly nonconvex) set $K \subset \U$.  (In particular, if $F = \nabla f$ for some $C^{(3)}$-smooth function $f \colon \U \to \R$, then the variational inequality reduces to the first-order condition for the optimization problem $\inf_K f$.)
Our aim is to apply Theorem \ref{active-set} (Active set method for generalized equations). 

To this end, we first suppose that the set $K$ is prox-regular at $-F(\bar u)$, and that 
$\bar u$ is in fact (following standard terminology \cite{fac_pang}) a {\em nondegenerate} solution:
\[
-F(\bar u) \in \mbox{ri}\,N_K(\bar u).
\]
(The practically important but much more complex case where nondegeneracy fails is analyzed in \cite{mirror}.)
Next, we assume that $K$ is $C^{(3)}$-{\em partly smooth} at $\bar u$ for $-F(\bar u)$ relative to some $C^{(3)}$-manifold $\M \subset K$ around $\bar u$, by which we mean that the normal cone mapping $N_K$ is inner semicontinuous at $\bar u$ relative to $\M$ and satisfies the {\em sharpness} condition
\[
\mbox{span}\,N_K(\bar u) = N_\M(\bar u).
\] 
In this case we have the local property \cite{ident}
\bmye \label{local-bundle}
\mbox{gph}\,N_K = \mbox{gph}\,N_\M \quad \mbox{around}~ \big(\bar u, -F(\bar u)\big).
\emye
(The graph of $N_\M$ is sometimes called the {\em normal bundle} for the manifold $\M$.)  It is precisely this property that underlies active set approaches to the original problem.

\subsubsection*{Example:  the basic projection algorithm}
To illustrate, suppose for this example that the set $K$ is convex. For any constant $a > 0$, by defining a continuous map $y_a \colon \U \to Q$ by
\[
y_a(u) = \mbox{Proj}_K \big(u - \frac{1}{a}F(u)\big)
\]
for points $u \in \U$, following a standard route \cite[Section 12.1]{fac_pang} we can convert the variational inequality into a fixed-point problem: $u$ solves $\mbox{VI}(K,F)$ if and only if $y_a(u) = u$. 
Indeed, the map above is just a special case of the map (\ref{restoration-Q}) that we defined for more general maximal monotone operators.  In the example of convex optimization, where $F$ is the gradient of a smooth convex function, the iteration
$u \leftarrow y_a(u)$ is just the projected gradient method.  Again specializing our earlier discussion,
by definition we have
\[
a\big( u - y_a(u) \big) - F(u)
~\in~
N_K\big(y_a(u)\big)
\]
for all points $u \in U$.  Under certain general conditions [Facchinei-Pang Thm 12.1.2], any sequence of iterates $u_k \in \U$ (for $k=0,1,2,\ldots$) must converge to the solution $\bar u$.  By continuity we have 
\[
y_k = y_a(u_k) \to \bar u,
\]
and also
\[
v_k ~=~ a( u_k - y_k ) - F(u_k) ~\to~ -F(\bar u),
\]
with
\[
v_k \in N_K(y_k).
\]
Hence, by the reduction property (\ref{local-bundle}), we deduce $y_k \in \M$ eventually:  in other words, the sequence $y_k$ {\em identifies} the active manifold $\M$.  In the case of the projected gradient method for convex optimization, this is exactly the result of \cite{Wright}, based on earlier work in \cite{Dunn87}.  The monograph \cite{fac_pang} includes a wide-ranging discussion of general solution algorithms for variational inequalities:  the identification behavior above holds for quite broad classes.  Global techniques often aim to drive down the natural {\em residual} $\|u - y_a(u)\|$ indirectly by driving down a merit function such as the {\em regularized gap}
\[
\theta_a(u) ~=~ \sup_{y \in K} \big\{ \ip{F(u)}{u-y} - \frac{a}{2}\| u-y \|^2 \big\}
\]
or the {\em D-gap} given by the function $\theta_a - \theta_b$, for some constant $b$ in the interval $(0,a)$ (see [Facchinei-Pang Propositions 10.3.7 and 10.3.8]).  Such approaches also lead to identification: if the iterates $u_k$ converge to some limit $\bar u$, and the corresponding residuals converge to zero, then continuity ensures that $\bar u$ is a fixed point of the map $y_a$ and hence a solution of the variational inequality, so under the partial smoothness assumptions above, the nearby sequence $y_a(u_k)$ eventually lies in the manifold $\M$.
\bigskip

Returning to the general case, property (\ref{local-bundle}) shows that, locally, around the solution 
$\bar u$, the original variational inequality $\mbox{VI}(K,F)$ is equivalent to the variational inequality $\mbox{VI}(\M,F)$:
\[
0 \in F(u) + N_\M(u).
\]
As we remarked, this reduction captures the essence of the active set philosophy:  we can henceforth focus on applying Theorem \ref{active-set} for the mapping $\Phi = F + N_\M$. 

Following \cite{lewis-liang}, the normal space mapping $N_\M$ is $C^{(2)}$-partly smooth at $\bar u$ for 
$-F(\bar u)$, with active manifold $\M$, so by an easy calculus rule \cite{lewis-liang}, the mapping 
$F + N_\M$ is $C^{(2)}$-partly smooth at $\bar u$ for $0$, as we need in order to apply 
Theorem~\ref{active-set}.  The theorem involves two more key ingredients.  The first is the transversality condition to check, which we can write
\[
 -\big(DF(\bar u)\big)^* y \in D^* N_\M \big(\bar u | -F(\bar u)\big)(y)
 \quad \Rightarrow \quad y=0,
\]
using routine coderivative calculus \cite[Example 10.43]{VA}.  The second is the linearized equation: given points $u \in \U$ near 
$\bar u$ and small $v \in F(u) + N_\M(u)$, we solve
\[
-v ~\in~ DF(u)(w-u) + DN_\M \big(u | v - F(u)\big)(w-u)
\]
for the new iterate $w$.  Here, once again we have applied some standard calculus, this time for derivatives \cite[Example 10.43]{VA}.  

In this light, we therefore next discuss how to compute the derivative and coderivative maps
\bmye \label{derivative-coderivative}
DN_\M(u | z) \quad \mbox{and} \quad D^*N_\M(u | z)
\emye
for $u \in \M$ and $z \in N_\M(u)$.  These are second-order concepts, since $N_\M = \partial \delta_\M$.  Indeed, using the notation of the Mordukhovich {\em generalized Hessian} \cite{boris-hessian}, we have 
\[
D^*N_\M = \partial^2\delta_\M.  
\]
The Hessian of a $C^{(2)}$-smooth function $f \colon \U \to \R$ has an important symmetry property that we can express as the relationship
\[
D(\nabla f) = D^*(\nabla f).
\]
By contrast, for general sets $\M$, there is no relationship between the derivative and coderivative mappings (\ref{derivative-coderivative}).  On the other hand, for fully amenable sets (and hence, in particular, for manifolds) we always have the {\em derivative-coderivative inclusion}
\[
DN_\M(u | z) \subset D^*N_\M(u | z);
\]
see \cite[Theorem 13.57]{VA}.  Our next result shows that, for manifolds, this inclusion in fact holds with equality:  informally, second derivatives of manifold indicator functions are symmetric.

\begin{thm}[Derivative-coderivative equality for manifolds] \label{equality} \mbox{} \\
Consider a
$C^{(2)}$-smooth manifold $\M \subset \U$ around a point 
$u \in \M$, and a normal vector \mbox{$v \in N_\M(u)$}.  Then the normal mapping 
$N_\M \colon \U \tto \U$ satisfies the derivative-coderivative equality
\[
DN_\M(u | v) = D^*N_\M(u | v).
\]
\end{thm}

\pf
By definition, we can suppose the manifold $\M$ to be defined as the set of common zeros near the point 
$u$ of \mbox{$C^{(2)}$-smooth} functions $h_i \colon \U \to \R$ (for $i=1,2,3,\ldots,m$), where the  gradients $\nabla h_i(u)$ are linearly independent, and that 
$v = \sum_i \lambda_i \nabla h_i(u)$ for some vector $\lambda \in \R^m$.  Define a corresponding self-adjoint linear map $H \colon \U \to \U$ by $H = \sum_i \lambda_i \nabla^2 h_i(u)$.  
By definition, for vectors $u,w \in \U$ we have
\[
z \in DN_\M(u | v)(w)
\]
if and only if 
\[
(w,z) \in T_{\mbox{\scriptsize gph}\,N_\M}(u,v),
\]
which is equivalent (by \cite[Theorem 2.8]{shanshan}) to
\bmye \label{tangent}
w \in T_\M(u) \quad \mbox{and} \quad z - Hw \in N_\M(u).
\emye
On the other hand, we have
\[
z \in D^*N_\M(u | v)(w)
\]
if and only if
\[
(z,-w) \in N_{\mbox{\scriptsize gph}\,N_\M}(u,v),
\]
which is also equivalent (by \cite[Theorem 2.8]{shanshan}) to condition (\ref{tangent}). 
\finpf

As an aside, we note an elegant consequence of this result:  again informally speaking, partly smooth {\em functions} have symmetric second derivatives.  To be precise, we gather together in the next result various properties of functions $f \colon \U \to \bR$ that are partly smooth in the sense of \cite{ident-arxiv}, the functional version of the definition for sets in \cite{ident}.  This idea is a slight variant of the original notion from \cite{Lewis-active}:  it corresponds closely to our current notion of partial smoothness applied to the subdifferential mapping $\partial f$.

\begin{cor}
Consider a $C^{(r)}$-smooth manifold $\M \subset \U$ around a point 
$\bar u \in \M$, for some degree $r>1$.  Suppose a subgradient $v$ of a function $f \colon \U \to \bR$ at $\bar u$ satisfies the nondegeneracy condition
\[
\bar v \in \mbox{\rm ri}\,\partial f(\bar u),
\]
as well as the following conditions:
\begin{itemize}
\item
$f$ is prox-regular and subdifferentially continuous at $\bar u$ for $\bar v$;
\item
$f$ agrees with some $C^{(r)}$-smooth function $\bar f \colon \U \to \R$ on $\M$ around $\bar u$;
\item
The affine span of the regular subdifferential $\hat \partial f(\bar u)$ is a translate of the normal space $N_\M(\bar u)$;
\item
For some neighborhood $W$ of $\bar v$, the mapping $u \mapsto \partial f(u) \cap W$ is inner semicontinuous at $\bar u$ relative to $\M$.
\end{itemize}
Then the subdifferential mapping $\partial f$ is $C^{(r-1)}$-partly smooth at $\bar u$ for $\bar v$, with active manifold $\M$, and, locally, satisfies
\[
\mbox{\rm gph}\, \partial f ~=~ \{(u,\nabla\bar f(u) + v) : u \in \M,~ v \in N_\M(u) \}
\quad \mbox{around}~ (\bar u,\bar v).
\]
Furthermore, the following derivative-coderivative equality holds:
\[
D(\partial f)(\bar u | \bar v) ~=~ D^*(\partial f)(\bar u | \bar v) ~=~ \partial^2 f(\bar u | \bar v)
 ~=~ \nabla^2 \bar f(\bar u) + D(N_\M)\big(\bar u | \bar v - \nabla f(\bar u)\big).
\]
\end{cor}

\noindent
{\bf Note:}  The bulleted properties (without the assumption of subdifferential continuity) together constitute the definition of {\em partial smoothness} of $f$ at $\bar u$ for $\bar v$.
\medskip

\pf
All except the final claim is simply a restatement from \cite{lewis-liang}.  The final claim follows from the fact that the graphs of the mappings $\partial f$ and $\nabla\bar f + N_\M$ agree around the point 
$(\bar u,\bar v)$, and hence these mappings have the same derivative and coderivative mappings.  The result now follows from derivative and coderivative calculus and Theorem \ref{equality}.
\finpf

Having calculated the relevant derivatives and coderivatives in Theorem \ref{equality}, we can now interpret the active set method described in Theorem \ref{active-set} for the variational inequality 
$\mbox{VI}(K,F)$.  Under the partial smoothness assumptions at the beginning of this section, the problem reduces locally around the solution $\bar u$ to variational inequality
\[
0 \in F(u) + N_\M(u),
\]
or equivalently, using the notation of the proof of Theorem \ref{equality}, to the system
\begin{eqnarray}
F(u) + \sum_{i=1}^m \lambda_i \nabla h_i(u) &=& 0 \label{first} \\
h_i(u) &=& 0 \quad (i=1,2,3,\ldots,m). \label{second}
\end{eqnarray}
in the variables $u \in \U$ and $\lambda \in \R^m$.  By assumption, there exists a unique vector 
$\bar\lambda \in \R^m$ such that the function $G = F + \sum_i \bar\lambda_i \nabla h_i$ is zero at 
$\bar u$.

As we have seen, the transversality condition is
\[
 -\big(DF(\bar u)\big)^* y \in D^* N_\M \big(\bar u | -F(\bar u)\big)(y)
 \quad \Rightarrow \quad y=0,
\]
or equivalently, by property (\ref{tangent}),
\[
y \in T_\M(\bar u) \quad \mbox{and} \quad 
DG(\bar u)^*y \in N_\M(\bar u)  \quad \Rightarrow \quad y=0.
\]
Not surprisingly, this simply amounts to the invertibility of the system (\ref{first}) and (\ref{second}) linearized around the solution $(\bar u,\bar  \lambda)$.  In particular, under these conditions the fact that $\bar u$ is an isolated solution is an automatic consequence.  In the case of optimization (where the map $F$ is a gradient), the condition reduces to nonsingularity of the Hessian of the Lagrangian function, projected onto the tangent space $T_\M(\bar u)$.

Turning to the linearized equation, we consider a current point $u \in \U$ near 
$\bar u$ and small vector $v \in F(u) + N_\M(u)$, which we can write as $L(u)$, where
\[
L = F + \sum_i \lambda_i \nabla h_i,  
\]
is the {\em Lagrangian}, for some unique multiplier vector $\lambda \in \R^m$.
We must solve
\[
-v ~\in~ DF(u)s + DN_\M \Big(u \Big| \sum_i \lambda_i \nabla h_i(u) \Big)s
\]
for the Newton step $s$, and then update $u \leftarrow u+s$.  We can write this equation equivalently, by property (\ref{tangent}), as the linear system
\[
s \in T_\M(u) \quad \mbox{and} \quad L(u) + DL(u)s \in N_\M(u),
\]
so this step amounts to seeking a critical point of the Lagrangian restricted to the tangent space, a familiar operation in sequential quadratic programming approaches.  

We can now interpret Theorem \ref{active-set}, after using this active-set approach to reduce to the concrete system (\ref{first}) and (\ref{second}).  We linearize this system around the current iterate 
$u$ and current Lagrange multiplier estimate $\lambda \in \R^m$, and solve for the new pair 
$(u',\lambda')$.  We then resort to exact feasibility restoration --- often unrealistic in practice, but nonetheless useful conceptually.  To be precise, we map the first component $u'$ back onto a new iterate $u$ on the active manifold $\M$ defined by the equations (\ref{second}):  conceptually, we could simply consider the nearest-point projection, especially if the original set $K$ is convex, since then the projections $\mbox{Proj}_K$ and $\mbox{Proj}_\M$ agree near the solution.  Finally, we estimate a new multiplier vector $\lambda \in \R^m$, for example by minimizing $\|F(u) + \sum_i \lambda_i \nabla h_i(u)\|$.  Repeating this process generates iterates $(u,\lambda)$ converging quadratically to the solution pair $(\bar u,\bar \lambda)$.

\subsubsection*{Acknowledgements}
Thanks to Jiayi Guo, who contributed to early versions of some of the ideas presented here.


\begin{thebibliography}{10}

\bibitem{Al-Khayyal-Kyparisis91}
F.~Al-Khayyal and J.~Kyparisis.
\newblock Finite convergence of algorithms for nonlinear programs and
  variational inequalities.
\newblock {\em J. Optim. Theory Appl.}, 70(2):319--332, 1991.

\bibitem{Burke90}
J.V. Burke.
\newblock On the identification of active constraints. {II}.\ {T}he nonconvex
  case.
\newblock {\em SIAM J. Numer. Anal.}, 27(4):1081--1103, 1990.

\bibitem{Burke-More88}
J.V. Burke and J.J. Mor{\'e}.
\newblock On the identification of active constraints.
\newblock {\em SIAM J. Numer. Anal.}, 25(5):1197--1211, 1988.

\bibitem{Calamai-More87}
P.H. Calamai and J.J. Mor{\'e}.
\newblock Projected gradient methods for linearly constrained problems.
\newblock {\em Math. Program.}, 39(1):93--116, 1987.

\bibitem{ident}
D.~Drusvyatskiy and A.S. Lewis.
\newblock Optimality, identifiability, and sensitivity.
\newblock {\em Math. Program.}, 147:467--498, 2014.

\bibitem{ident-arxiv}
D.~Drusvyatskiy and A.S. Lewis.
\newblock Optimality, identifiability, and sensitivity.
\newblock \\ arXiv:1207.6628, 2014.

\bibitem{Dunn87}
J.C. Dunn.
\newblock On the convergence of projected gradient processes to singular
  critical points.
\newblock {\em J. Optim. Theory Appl.}, 55(2):203--216, 1987.

\bibitem{fac_pang}
F.~Facchinei and J.~Pang.
\newblock {\em Finite {D}imensional {V}ariational {I}nequalities and
  {C}omplementarity {P}roblems}.
\newblock Springer Series in Operations Research, Springer-Verlag, New York,
  2003.

\bibitem{mirror}
J.~Fadili, J.~Malick, and G.~Peyr\'{e}.
\newblock Sensitivity analysis for mirror-stratifiable convex functions.
\newblock {\em SIAM J. Optim.}, 28:2975--3000, 2018.

\bibitem{Ferris91}
M.C. Ferris.
\newblock Finite termination of the proximal point algorithm.
\newblock {\em Math. Program. Ser. A}, 50(3):359--366, 1991.

\bibitem{Flam92}
S.D. Fl{\aa}m.
\newblock On finite convergence and constraint identification of subgradient
  projection methods.
\newblock {\em Math. Program.}, 57:427--437, 1992.

\bibitem{HS09}
W.~Hare and C.~Sagastiz{\'a}bal.
\newblock Computing proximal points of nonconvex functions.
\newblock {\em Math. Program.}, 116(1-2, Ser. B):221--258, 2009.

\bibitem{Lee}
J.M. Lee.
\newblock {\em Introduction to Smooth Manifolds}.
\newblock Springer, New York, 2003.

\bibitem{dual_av}
S.~Lee and S.J. Wright.
\newblock Manifold identification in dual averaging for regularized stochastic
  online learning.
\newblock {\em J. Mach. Learn. Res.}, 13:1705--1744, 2012.

\bibitem{LOS}
C.~Lemar\'{e}chal, F.~Oustry, and C.~Sagastiz\'{a}bal.
\newblock The {U}-lagrangian of a convex function.
\newblock {\em Transactions of the American Mathematical Society},
  352:711--729, 2000.

\bibitem{LS97}
C.~Lemar\'{e}chal and C.~Sagastiz\'{a}bal.
\newblock Practical aspects of the {M}oreau-{Y}osida regularization:
  theoretical preliminaries.
\newblock {\em SIAM J. Optim.}, 7:367---385, 1997.

\bibitem{Lewis-active}
A.S. Lewis.
\newblock Active sets, nonsmoothness, and sensitivity.
\newblock {\em SIAM J. Optim.}, 13:702--725, 2002.

\bibitem{lewis-liang}
A.S. Lewis and Jingwei Liang.
\newblock Partial smoothness and constant rank.
\newblock \\ {\tt arXiv:1807.03134}, 2018.

\bibitem{prx_lin}
A.S. Lewis and S.J. Wright.
\newblock A proximal method for composite minimization.
\newblock {\em Math. Program.}, 158:501--546, 2016.

\bibitem{shanshan}
A.S. Lewis and S.~Zhang.
\newblock Partial smoothness, tilt stability, and generalized {H}essians.
\newblock {\em SIAM J. Optim.}, 23:74--94, 2013.

\bibitem{liang-phd}
J.~Liang.
\newblock {\em Convergence Rates of First-Order Operator Splitting Methods}.
\newblock PhD thesis, University of Caen Normandie, 2016.

\bibitem{fb_ps}
J.~Liang, J.~Fadili, and G.~Peyr\'{e}.
\newblock Local linear convergence of forward--backward under partial
  smoothness.
\newblock In Z.~Ghahramani, M.~Welling, C.~Cortes, N.D. Lawrence, and K.Q.
  Weinberger, editors, {\em Advances in Neural Information Processing Systems
  27}, pages 1970--1978. Curran Associates, Inc., 2014.

\bibitem{fadili-activity}
J.~Liang, J.~Fadili, and G.~Peyr\'{e}.
\newblock Activity identification and local linear convergence of
  forward-backward-type methods.
\newblock {\em SIAM J. Optim.}, 27:408--437, 2017.

\bibitem{liang}
J.~Liang, J.~Fadili, and G.~Peyr\'{e}.
\newblock Local linear convergence analysis of primal-dual splitting methods.
\newblock {\em Optimization}, 67:821--853, 2018.

\bibitem{activ_ident}
J.~Liang, J.~Fadili, G.~Peyr\'{e}, and R.~Luke.
\newblock Activity identification and local linear convergence of
  {D}ouglas-{R}achford/{ADMM} under partial smoothness.
\newblock In J.-F. Aujol, M.~Nikolova, and N.~Papadakis, editors, {\em Scale
  Space and Variational Methods in Computer Vision}, volume 9087 of {\em
  Lecture Notes in Computer Science}, pages 642--653. Springer International
  Publishing, 2015.

\bibitem{MC02}
R.~Mifflin and C.~Sagastiz{\'a}bal.
\newblock Proximal points are on the fast track.
\newblock {\em Journal of Convex Analysis}, 9:563---579, 2002.

\bibitem{MC03}
R.~Mifflin and C.~Sagastiz{\'a}bal.
\newblock Primal-dual gradient structured functions: second-order results;
  links to epi-derivatives and partly smooth functions.
\newblock {\em SIAM J. Optim.}, 13:1174--1194, 2003.

\bibitem{MC04}
R.~Mifflin and C.~Sagastiz{\'a}bal.
\newblock {$\cal V\cal U$}-smoothness and proximal point results for some
  nonconvex functions.
\newblock {\em Optim. Methods Softw.}, 19:463--478, 2004.

\bibitem{MC05}
R.~Mifflin and C.~Sagastiz{\'a}bal.
\newblock A {$\cal{VU}$}-algorithm for convex minimization.
\newblock {\em Math. Program.}, 104:583--608, 2005.

\bibitem{miller-malick}
S.A. Miller and J.~Malick.
\newblock Newton methods for nonsmooth convex minimization: connections among
  {$\mathcal U$}-{L}agrangian, {R}iemannian {N}ewton and {SQP} methods.
\newblock {\em Math. Program.}, 104:609--633, 2005.

\bibitem{boris-hessian}
B.S. Mordukhovich.
\newblock Sensitivity analysis in nonsmooth optimization.
\newblock In {\em Theoretical {A}spects of {I}ndustrial {D}esign
  ({W}right-{P}atterson {A}ir {F}orce {B}ase, {OH}, 1990)}, pages 32--46. SIAM,
  Philadelphia, PA, 1992.

\bibitem{Mord_1}
B.S. Mordukhovich.
\newblock {\em Variational Analysis and Generalized Differentiation I: Basic
  Theory}.
\newblock Grundlehren der mathematischen Wissenschaften, Vol 330, Springer,
  Berlin, 2006.

\bibitem{nocedal_wright}
J.~Nocedal and S.J. Wright.
\newblock {\em Numerical Optimization}.
\newblock Springer Series in Operations Research and Financial Engineering.
  Springer, New York, second edition, 2006.

\bibitem{VA}
R.T. Rockafellar and R.J-B. Wets.
\newblock {\em Variational {A}nalysis}.
\newblock Grundlehren der mathematischen Wissenschaften, Vol 317, Springer,
  Berlin, 1998.

\bibitem{claudia-icm}
C.~Sagastiz\'{a}bal.
\newblock A {$\cal{VU}$}-point of view of nonsmooth optimization.
\newblock In {\em Proceedings of the International Congress of Mathematicians,
  Rio de Janeiro}, volume~3, pages 3785--3806, 2018.

\bibitem{Vaiter2017}
S.~Vaiter, C.~Deledalle, J.~Fadili, G.~Peyr\'{e}, and C.~Dossal.
\newblock The degrees of freedom of partly smooth regularizers.
\newblock {\em Annals of the Institute of Statistical Mathematics},
  69:791--832, 2017.

\bibitem{fadili-model}
S.~Vaiter, G.~Peyr\'{e}, and J.~Fadili.
\newblock Model consistency of partly smooth regularizers.
\newblock {\em IEEE Trans. Inform. Theory}, 64:1725--1737, 2018.

\bibitem{Wright}
S.J. Wright.
\newblock Identifiable surfaces in constrained optimization.
\newblock {\em SIAM J. Control Optim.}, 31:1063--1079, 1993.

\end{thebibliography}

\def\cprime{$'$} \def\cprime{$'$}

\end{document}